\documentclass{article}

    \PassOptionsToPackage{numbers, compress}{natbib}
\usepackage[preprint]{neurips_2019}

\usepackage[utf8]{inputenc} % allow utf-8 input
\usepackage[T1]{fontenc}    % use 8-bit T1 fonts
\usepackage{hyperref}       % hyperlinks
\usepackage{url}            % simple URL typesetting
\usepackage{booktabs}       % professional-quality tables
\usepackage{amsfonts}       % blackboard math symbols
\usepackage{nicefrac}       % compact symbols for 1/2, etc.
\usepackage{microtype}      % microtypography
\usepackage{amssymb,amsmath, amsthm}
\usepackage{algorithm,algorithmic}
\usepackage{color}
\usepackage{graphicx}
\usepackage{subfigure}
\newtheorem{theorem}{Theorem}
\newtheorem{remark}{Remark}
\newtheorem{assumption}{Assumption}

\title{Projected Stein Variational Newton: A Fast and Scalable Bayesian Inference Method \\ in High Dimensions}

\author{%
	Peng~Chen, Keyi~Wu, Joshua~Chen, Thomas~O’Leary-Roseberry, Omar~Ghattas\\
	Oden Institute for Computational Engineering and Sciences \\
	The University of Texas at Austin \\
	Austin, TX 78712.\\
	\texttt{\{peng, keyi, joshua, tom, omar\}@oden.utexas.edu} \\
}

\newcommand{\bE}{{\mathbb{E}}}
\newcommand{\bH}{{\mathbb{H}}}
\newcommand{\bN}{{\mathbb{N}}}
\newcommand{\bR}{{\mathbb{R}}}

\newcommand{\cD}{\mathcal{D}}

\newcommand{\cH}{\mathcal{H}}
\newcommand{\cJ}{\mathcal{J}}
\newcommand{\cN}{\mathcal{N}}
\newcommand{\cT}{\mathcal{T}}

\newcommand{\bsc}{\boldsymbol{c}}
\newcommand{\bsg}{\boldsymbol{g}}

\newcommand{\pc}[1]{{{\color{red}{#1}}}}

\newcommand{\beq}{\begin{equation}}
\newcommand{\eeq}{\end{equation}}

\begin{document}
	
	\maketitle
	\begin{abstract}
		We propose a projected Stein variational Newton (pSVN) method for high-dimensional Bayesian inference. To address the curse of dimensionality, we exploit the intrinsic low-dimensional geometric structure of the posterior distribution in the high-dimensional parameter space via its Hessian (of the log posterior) operator and perform a parallel update of the parameter samples projected into a low-dimensional subspace by an SVN method. The subspace is adaptively constructed using the eigenvectors of the averaged Hessian at the current samples. We demonstrate fast convergence of the proposed method, complexity independent of the parameter and sample dimensions, and parallel scalability.
		%% OG: shouldn't use term "scalability" to refer to arithmetic complexity, since that caused confusion last time. 
		
	\end{abstract}
	
	\section{Introduction}
	\label{introduction}
	Bayesian inference provides an optimal probability formulation for learning complex models from observational or experimental data under uncertainty by updating the model parameters from their prior distribution to a posterior
	%% OG: need to say "parameters" rather than parameter, since you're talking to a finite dimensional crowd. They won't understand that parameter is an (infinite dim) parameter field. 
	distribution \cite{Stuart10}. In Bayesian inference we typically face the task of drawing samples from the posterior probability distribution to compute various statistics of some given quantities of interest. 
	However, this is often prohibitive when the posterior distribution is high-dimensional; many conventional methods for Bayesian inference suffer from the curse of dimensionality, i.e., computational complexity grows exponentially or convergence deteriorates with increasing parameter dimension.
	
	To address this curse of dimensionality, several efficient and dimension-independent methods have been developed that exploit the intrinsic properties of the posterior distribution, such as its smoothness, sparsity, and intrinsic low-dimensionality. Markov chain Monte Carlo (MCMC) methods exploiting geometry of the log-likelihood function have been developed \cite{GirolamiCalderhead11, MartinWilcoxBursteddeEtAl12, StadGhatt2014, CuiLawMarzouk16, BESKOS2017327}, providing more effective sampling than black-box MCMC. For example, the DILI MCMC method \cite{CuiLawMarzouk16} uses the low rank structure of the Hessian of the negative log likelihood in conjunction with operator-weighted proposals that are well-defined on function space to yield a sampler whose performance is dimension-independent and effective at capturing information provided by the data. However, despite these enhancements, MCMC methods remain prohibitive for problems with expensive-to-evaluate likelihoods (i.e., involving complex models) and in high parameter dimensions. Deterministic sparse quadratures were developed in \cite{SchwabStuart12, SchillingsSchwab2013, ChenSchwab2015} and shown to converge rapidly with dimension-independent rates for smooth and sparse problems. However, the fast convergence is lost when the posterior has significant local variations, despite enhancements with Hessian-based transformation  \cite{SchillingsSchwab16, ChenVillaGhattas2017}. 
	
	Variational inference methods reformulate the sampling problem as an optimization problem that approximates the posterior by minimizing its Kullback--Leibler divergence with a transformed prior \cite{MarzoukMoselhyParnoEtAl16, LiuWang16, BleiKucukelbirMcAuliffe17}, which can be potentially much faster than MCMC. In particular, Stein variational methods, which seek a composition of a sequence of simple transport maps represented by kernel functions using gradient descent (SVGD) \cite{LiuWang16, ChenMackeyGorhamEtAl18, LiuZhu18} and especially Newton (SVN) \cite{DetommasoCuiMarzoukEtAl18} optimization methods, are shown to achieve fast convergence in relatively low dimensions.
	However, these variational optimization methods can again become deteriorated in convergence and accuracy in high dimensions. The curse of dimensionality can be partially addressed by a localized SVGD on Markov blankets, which relies on a conditional independence structure of the target distribution \cite{ZhuoLiuShiEtAl17, WangZengLiu18}. 
	
	\textbf{Contributions}: In this work, we develop a projected Stein variational Newton method (pSVN) to tackle the challenge of high-dimensional Bayesian inference by exploiting the intrinsic low-dimensional geometric structure of the posterior distribution (where it departs from the prior), as characterized by the dominant spectrum of the prior-preconditioned Hessian of the negative log likelihood. This low-rank structure, or fast decay of eigenvalues of the preconditioned Hessian, has been proven  for some inference problems and commonly observed in many others with complex models \cite{Bui-ThanhGhattas2012, Bui-ThanhGhattasMartinEtAl2013, SpantiniSolonenCuiEtAl15, IsaacPetraStadlerEtAl15, CuiLawMarzouk16, ChenVillaGhattas2017, ChenVillaGhattas18, BashirWillcoxGhattasEtAl08, chenGhattas18b}. By projecting the parameters into this data-informed low-dimensional subspace and applying the SVN in this subspace, we can effectively mitigate the curse of dimensionality.
	We demonstrate fast convergence of pSVN that is independent of the number of parameters and samples. In particular, in two (both linear and nonlinear) experiments we show that the intrinsic dimension is a few (6) and a few tens (40) with the nominal dimension over 1K and 16K, respectively. We present a scalable parallel implementation of pSVN that yields rapid convergence, minimal communication, and low memory footprint, thanks to this low-dimensional projection. 
	
	Below, we present background on Bayesian inference and Stein variational methods in Section \ref{sec:background}, develop the projected Stein variational Newton method in Section \ref{sec:pSVN}, and provide numerical experiments in Section \ref{sec:numerics}.
	
	% \section{Stein variational Bayesian inference}
	\section{Background}
	\label{sec:background}
	% We overview the formulation of Bayesian inference and Stein variational methods in this section.
	\subsection{Bayesian inference}
	
	We consider a random parameter $x \in \bR^d$, $d \in \bN$, with a prior probability density function $p_0:\bR^d \to \bR$, and noisy observational data $y$ of a parameter-to-observable map $f:\bR^d \to \bR^s$, $s\in \bN$, i.e., 
	\beq
	y = f(x) + \xi, 
	\eeq
	where $\xi \in \bR^s$ represents observation noise with probability density function $p_\xi: \bR^s \to \bR$.
	%we consider Gaussian noise $\xi \sim \cN(0, \Gamma)$ with s.p.d. covariance $\Gamma \in \bR^{s\times s}$. 
	% For simplicity, we also consider that the prior is Gaussian with distribution $\cN(\bar{x}, \Gamma_0)$ with mean $\bar{x}$ and s.p.d. covariance $\Gamma_0$. 
	The posterior density $p(\cdot|y):\bR^d \to \bR$ of $x$ conditioned on the data $y$ is given by Bayes' rule  
	\beq\label{eq:Bayes_x}
	p(x|y) = \frac{1}{Z} p_y(x), \quad \text{ where } p_y(x) := p_\xi(y - f(x)) \, p_0(x),
	\eeq 
	and the normalization constant $Z$, typically $Z\neq 1$ if $p_\xi$ or $p_0$ is known up to a constant, is given by 
	\beq\label{eq:potential}
	Z := \bE_{p_0}[p_\xi(y - f(x))] = \int_{\bR^d} p_\xi(y-f(x)) p_0(x)dx.
	\eeq
	%Note that $||y||_\Gamma^2 = y^T \Gamma^{-1} y$. %\textcolor{red}{OG: why not use the standard norm symbol, $\parallel$?} 
	In practice, $Z$ is computationally intractable, especially for large $d$. 
	%One seeks to draw samples from the posterior, whose probability density, known up to the constant $Z$, is denoted by $p_y=e^{-\eta_y}p_0$.
	
	\subsection{Stein variational methods}
	\label{sec:SVM}
	While sampling from the prior is tractable, sampling from the posterior is a great challenge. One method to sample from the posterior is to find a transport map $T: \bR^d \to \bR^d$ in a certain function class $\cT$ that pushes forward the prior to the posterior by minimizing a Kullback--Leibler (KL) divergence 
	\beq\label{eq:DKL}
	\min_{T\in \cT} \mathcal{D}_{\text{KL}} (T_* p_0 | p_y).
	\eeq
	% where $T_*$ denotes the pushforward map. with its dual pullback map denoted as $T^*$, the KL divergence is given by
	% \beq
	% \mathcal{D}_{\text{KL}} (T_* p_0 | p_y)= \mathcal{D}_{\text{KL}} (p_0 | T^* p_y) = \int_{\bR^d} \log\left(\frac{p_0(x)}{p_y(T(x))\text{det}(\nabla T(x))}\right)p_0(x)dx,
	% \eeq
	% where $T_* \mu_0$ represents the pushforward measure such that 
	% \beq
	% \int_{\bR^d} g(x) dT_* \mu_0(x) = \int_{\bR^d} g \circ T(x) d\mu_0(x)
	% \eeq
	% for any $T_* \mu_0$-measurable function $g$. 
	Stein variational methods \cite{LiuWang16, DetommasoCuiMarzoukEtAl18} simplify the minimization of \eqref{eq:DKL} for one possibly very complex and nonlinear transport map $T$ to a sequence of simpler transport maps that are perturbations of the identity, i.e., $T = T_L \circ T_{L-1} \circ \cdots \circ T_2 \circ T_1$, $L \in \bN$, where
	\beq\label{eq:mapl}
	T_l(x) = I(x) + \varepsilon Q_l(x), \quad l = 1, \dots, L,
	\eeq
	with $I(x) = x$, step size $\varepsilon$, and perturbation map $Q_l:\bR^d \to \bR^d$. Let $p_l$ denote the pushforward density $p_l := (T_l \circ \cdots \circ T_1)_* p_0$. For $l = 1, 2, \dots$, we define a cost functional $\cJ_l(Q)$ as
	\beq
	\cJ_l(Q) := \mathcal{D}_{\text{KL}} ((I + Q)_* p_{l-1} | p_y).
	\eeq
	Then at step $l$, Stein variational methods lead to 
	\beq
	Q_l = - \cH^{-1}_l \nabla \cJ_l(0),
	\eeq
	where $\nabla \cJ_l(0) : \bR^d \to \bR^d$ is the Fr\'echet derivative of $\cJ_l(Q)$ evaluated at $Q = 0$, and $\cH_l$ is a preconditioner. For the SVGD method \cite{LiuWang16}, $\cH_l = I$, while for the SVN method \cite{DetommasoCuiMarzoukEtAl18}, $\cH_l \approx \nabla^2 \cJ_l(0)$, an approximation of the Hessian of the cost functional $\nabla^2 \cJ_l(0)$.
	
	Given basis functions $k_n:\bR^d \to \bR$, $n =1, \dots, N$, an ansatz representation of $Q_l$ is defined as 
	\beq\label{eq:Galerkin}
	Q_l(x) = \sum_{n = 1}^N c_n k_n(x), 
	\eeq
	where $c_n \in \bR^d$, $n =1, \dots, N$, are unknown coefficient vectors. It is shown in \cite{DetommasoCuiMarzoukEtAl18} that the coefficient vector $\bsc = (c_1^{\top}, \dots, c_N^{\top})^{\top} \in \bR^{Nd}$ is a solution of the linear system
	\beq\label{eq:system}
	\bH \bsc = - \bsg,
	\eeq
	where $\bsg = (g_1^{\top}, \dots, g_N^{\top})^{\top} \in \bR^{Nd}$ is the gradient vector given by
	\beq\label{eq:gradient}
	g_{m} := \bE_{p_{l-1}} [-\nabla_x \log(p_y) k_m - \nabla_x k_m], \quad m = 1, \dots, N,
	\eeq
	and $\bH \in \bR^{Nd\times Nd}$ is the Hessian matrix, specified as the identity for SVGD \cite{LiuWang16}, which leads to $c_n = -g_n$, $n = 1, \dots, N$, while for SVN it is given with $mn$-block $\bH_{mn} \in \bR^{d\times d}$ by \cite{DetommasoCuiMarzoukEtAl18}
	\beq\label{eq:Hessian}
	\bH_{mn} :=  \bE_{p_{l-1}} [-\nabla^2_x \log(p_y) k_n k_m + \nabla_x k_n (\nabla_x k_m)^{\top}], \quad m, n = 1, \dots, N.
	\eeq
	% To solve the coupled $Nd \times Nd$ system \eqref{eq:system}, a diagonal approximation is used in \cite{DetommasoCuiMarzoukEtAl18}, i.e., 
	% \beq
	% \bH_{mm} c_m = - g_m, \quad m =1, \dots, N.
	% \eeq
	At each step $l = 1, 2, \dots$, the expectation $\bE_{p_{l-1}}[\cdot]$ in \eqref{eq:gradient} and \eqref{eq:Hessian} are approximated by the sample average approximation with samples $x^{l-1}_1, \dots, x^{l-1}_N$, which are drawn from the prior at $l = 1 $ and pushed forward by \eqref{eq:mapl} once the coefficients $c_1, \dots, c_N$ are obtained. We remark that in the original SVGD method \cite{LiuWang16}, the samples are moved with the simplified perturbation $Q_l(x_m) = c_m$.  
	
	In both \cite{LiuWang16} and \cite{DetommasoCuiMarzoukEtAl18}, the basis functions $k_n(x)$
	are specified by a suitable kernel function $k_n(x) = k(x, x')$ at $x' = x_n$, $n = 1, \dots, N$, e.g., a Gaussian kernel given by
	\beq\label{eq:kernel}
	k(x, x') = \exp\left(-\frac{1}{2} (x-x')^{\top} M (x-x')\right),
	\eeq
	where $M$ is a metric that measures the distance between $x$ and $x' \in \bR^d$. In \cite{LiuWang16}, it is specified as rescaled identity matrix $\alpha I$ for $\alpha > 0$ depending on the samples, while in \cite{DetommasoCuiMarzoukEtAl18}, $M$ is given by $
	M = \bE_{p_{l-1}} [ -\nabla^2_x \log(p_y)]/d$ to account for the geometry of the posterior by averaged Hessian information. This was shown to accelerate convergence for both SVGD and SVN compared to $\alpha I$. We remark that for high-dimensional complex models where a direct computation of the Hessian $\nabla^2_x \log(p_y)$ is not tractable, its low-rank decomposition by randomized algorithms can be applied.
	
	%We remark that in high dimensions, the samples $x_1, \dots, x_N$ are often far from each other, which leads to $k_n(x_m) \approx \triangle_{mn}$ and $P_l(x_m) \approx c_m$ by \eqref{eq:Galerkin}.

	\section{Projected Stein variational Newton}
	\label{sec:pSVN}
	% In this section, we present pSVN to address high-dimensional Bayesian inference problems.
	
	\subsection{Dimension reduction by projection}
	Stein variational methods suffer from the curse of dimensionality, i.e., the sample estimate (e.g., for variance) deteriorates  considerably in high dimensions because the global kernel function \eqref{eq:kernel} cannot represent the transport map well, as shown in \cite{ZhuoLiuShiEtAl17, WangZengLiu18} for SVGD.
	This challenge can be alleviated in moderate dimensions by a suitable choice of the metric in \eqref{eq:kernel} as demonstrated in \cite{DetommasoCuiMarzoukEtAl18}. However it is still present when the dimension becomes high. An effective method to tackle this difficulty, which relies on conditional independence of the posterior density, uses local kernel functions defined over a Markov blanket with much lower dimension, thus achieving effective dimension reduction \cite{ZhuoLiuShiEtAl17, WangZengLiu18}. 
	
	In many applications, even when the nominal dimension of the parameter is very high, the intrinsic parameter dimension informed by the data is typically low, i.e., the posterior density is effectively different from the prior density only in a low-dimensional subspace  \cite{Bui-ThanhGhattas2012, Bui-ThanhGhattasMartinEtAl2013, SpantiniSolonenCuiEtAl15, IsaacPetraStadlerEtAl15, CuiLawMarzouk16, ChenVillaGhattas2017, ChenVillaGhattas18, BashirWillcoxGhattasEtAl08}. This is because: (i) the prior $p_0$ may have correlation in different dimensions, 
	(ii) the parameter-to-observable map $f$ may be smoothing/regularizing, (iii) the data $y$ may not be very informative, or a combined effect. Let $\Psi = (\psi_1, \dots, \psi_r) \in \bR^{d\times r}$ denote the basis of a subspace of dimension $r \ll d$ in $\bR^d$. Then we can project the parameter $x$ with mean $\bar{x}$ into this subspace as
	\beq\label{eq:projection}
	x^r = \bar{x} + P_r (x-\bar{x})= \bar{x} + \sum_{i =1}^r \psi_i (\psi_i, (x-\bar{x}))_H = \bar{x} + \sum_{i =1}^r \psi_i w_i = \bar{x} + \Psi w,
	\eeq
	where $w = (w_1, \dots, w_r) \in \bR^r$ is a vector of coefficients $w_i = (\psi_i, x-\bar{x})_H$ of the projection of $x -\bar{x}$ to $\psi_i$ in a suitable norm $H$, e.g., $(\psi_i, x-\bar{x})_H = \psi_i^T \Gamma_0^{-1} (x-\bar{x})$ where $\Gamma_0$ is the prior covariance of $x$ and $\psi_i^T \Gamma_0^{-1} \psi_j = \delta_{ij}$.  We define the projected posterior as 
	\beq\label{eq:Bayes_w}
	p^r(x|y) = \frac{1}{Z^r} p_y^r(x), \; \text{ where } p_y^r(x) = p_\xi(y - f(x^r)) p_0(x) \text{ and } Z^r = \bE_{p_0}[p_\xi(y - f(x^r))].
	\eeq
	Then we can establish convergence under the following assumption. We define $||\cdot||_X$ as a suitable norm, e.g., $||x||_X^2 = x^T X x$ with $X = I$, the identity matrix or a mass matrix discretized from identity operator in finite dimension approximation space in our numerical experiments.
	\begin{assumption}\label{ass:f}
		For Gaussian noise $\xi \in \cN(0, \Gamma)$ with s.p.d.\ covariance $\Gamma \in \bR^{s\times s}$. Let $||v||_\Gamma := (v^T \Gamma^{-1} v)^{1/2}$ for any $v\in \bR^s$. Assume
		there exists a constant $C_f > 0$ such that for any $x^r$ in \eqref{eq:projection}
		% For every $a > 0$, there is a constant $C_a \in \bR$ such that ($||\cdot||_X$ is a suitable norm in $\bR^d$)
		\beq
		\bE_{p_0} [||f(x^r)||_\Gamma] \leq C_f \; \text{ and } \; \bE_{p_0} [||f(x)||_\Gamma] \leq C_f.
		\eeq
		%\exp(a||x||_X^2 + C_a).
		For every $b > 0$, assume there is $C_b>0$ such that for all $x_1, x_2$ with $\max\{||x_1||_X, ||x_2||_X\} < b$, 
		\beq
		||f(x_1) - f(x_2)||_\Gamma \leq C_b ||x_1 - x_2||_X.
		\eeq
		% \item For every $a > 0$, and $C_a > 0$, there exists constant $C_p$ such that 
		% \beq
		% \int_{\bR^d} \exp(a||x||_X^2 + C_a) p_0(x) dx \leq C_p.
		% \eeq
		% This light tail condition is satisfied for Gaussian prior by Fernique's theorem \cite{Stuart10}.
		
	\end{assumption}
	
	We state the convergence result for the projected posterior density in the following theorem, whose proof is presented in Appendix A.
	\begin{theorem}\label{thm:convergence}
		Under Assumption \ref{ass:f}, there exists a constant $C$ independent of $r$ such that 
		\beq
		\cD_{\text{KL}} (p(x|y)\, |\,p^r(x|y)) \leq C ||x - x^r||_X.
		\eeq
		\begin{remark}
			Theorem \ref{thm:convergence} indicates that the projected posterior converges to the full one as along as the projected parameter converges in $X$-norm, and that the convergence of the former is bounded by the latter. In practical applications, the former may converge faster than the latter because it only depends on the data-informed subspace while the latter is measured in data-independent $X$-norm.
		\end{remark}
	\end{theorem}
	
	\subsection{Projected Stein variational Newton}
	Let $p_0^r$ denote the prior densities for $x^r $ in \eqref{eq:projection}. Let $x^\perp = x - x^r$. Then the prior is decomposed as 
	\beq
	p_0(x) = p_0^r(x^r) p_0^\perp(x^\perp| x^r),
	\eeq
	where $p_0^\perp(x^\perp| x^r)$ is a conditional density, which becomes $p_0^\perp(x^\perp)$ if $p_0$ is a Gaussian density. Then the projected posterior density $p_y^r(x)$ in \eqref{eq:Bayes_w} becomes 
	\beq
	p_y^r(x) = p_\xi(y - f(x^r)) p_0^r(x^r) p_0^\perp(x^\perp|x^r),
	\eeq
	so that sampling from $p_y^r(x)$ can be realized by sampling from $p_y^r(x^r) = p_\xi(y - f(x^r)) p_0^r(x^r)$ for $x^r$ and from $p_0^\perp(x^\perp|x^r)$ for $x^\perp$ conditioned on $x^r$ (or from $p_0^\perp(x^\perp)$ if $p_0$ is Gaussian). To sample from the posterior, we can sample $x$ from the prior, decompose it as $x = x^r + x^\perp$, freeze $x^\perp$, push $x^r$ to $x_y^r$ as a sample from $p_y^r(x^r)$, and construct the posterior sample as $x_y = x_y^r + x^\perp$.
	
	To sample from $p_y^r(x^r)$ in the projection subspace, we seek a transport map $T$ that pushes forward $p_0^r(x^r)$ to $p_y^r(x^r)$ by minimizing the KL divergence between them. Since the randomness of $x^r = \bar{x} + \Psi w$ is fully represented by $w$ given the projection basis $\Psi$, we just need to find a transport map that pushes forward $\pi_0(w) = p_0^r(x^r)$ to $\pi_y(w) = p_y^r(x^r)$ in the (coefficient) parameter space $\bR^r$, where $r \ll d$. Similarly in the full space, we look for a composition of a sequence of maps $T = T_L \circ T_{L-1} \circ \cdots \circ T_2 \circ T_1$, $L \in \bN$, with
	\beq\label{eq:mapl_w}
	T_l(w) = I(w) + \varepsilon Q_l(w), \quad l = 1, \dots, L,
	\eeq
	where the perturbation map $Q_l$ is represented by the basis functions $k_n:\bR^r \to \bR$, $n = 1, \dots, N$, as 
	\beq\label{eq:Galerkin_w}
	Q_l(w) = \sum_{n=1}^N c_n k_n(w) ,
	\eeq
	Then the coefficient vector $\bsc = ((c_1)^{\top}, \dots, (c_N)^{\top})^{\top} \in \bR^{Nr}$ is the solution of the linear system 
	\beq\label{eq:system_w}
	\bH \bsc = - \bsg.
	\eeq
	Here the  $m$-th component of the gradient $\bsg$ is defined as
	\beq\label{eq:gradient_m}
	g_m := \bE_{\pi_{l-1}}[-\nabla_w \log(\pi_y) k_m - \nabla_w k_m],
	\eeq
	and the $mn$-th component of the Hessian $\bH$ for pSVN is defined as 
	\beq\label{eq:Hessian_mn}
	\bH_{mn} :=  \bE_{\pi_{l-1}} [-\nabla^2_w \log(\pi_y) k_n k_m + \nabla_w k_n (\nabla_w k_m)^{\top}]. 
	\eeq
	The expectations in \eqref{eq:gradient_m} and \eqref{eq:Hessian_mn} are evaluated by sample average approximation at samples $w^{l-1}_1, \dots, w_N^{l-1}$, which are drawn from $\pi_0$ for $l=1$ and pushed forward by \eqref{eq:mapl_w} as $w^l_n = T(w^{l-1}_n)$, $n = 1, \dots, N$.
	By the definition of the projection \eqref{eq:projection}, we have
	\beq\label{eq:gradient_w}
	\nabla_w \log(\pi_y(w)) = \Psi^{\top} \nabla_x \log(p^r_y(x^r)),
	\text{ and }
	\nabla^2_w \log(\pi_y(w)) = \Psi^{\top} \nabla_x^2 \log(p^r_y(x^r)) \Psi.
	\eeq
	For the basis functions $k_n$, $n = 1, \dots, N$, we use a Gaussian kernel $k_n(w) = k(w,w_n)$ defined as in \eqref{eq:kernel}, with the metric $M$ given by  
	% \beq\label{eq:kernel_w}
	% k_n(w) = \exp\left(-\frac{1}{2}(w-w_n)^{\top} M^w (w-w_n)\right), 
	% \eeq
	% where the metric $M^w$ is given by 
	an averaged Hessian at the current samples $w_1^{l-1}, \dots, w_N^{l-1}$, i.e., 
	\beq\label{eq:metric}
	M = -\frac{1}{r} \bE_{\pi_{l-1}} [\nabla_w^2 \log(\pi_y)] \approx -\frac{1}{rN} \sum_{n=1}^N \nabla^2_w \log(\pi_y(w_n^{l-1})).
	\eeq
	We remark that the projected system \eqref{eq:system_w} is of size $Nr \times Nr$, which is a considerable reduction from the full system \eqref{eq:system} of size $Nd \times Nd$, since $r \ll d$. To further reduce the size of the coupled system \eqref{eq:system_w}, we use a classical ``mass-lumping'' technique to decouple it as $N$ systems of size $r \times r$ 
	\beq\label{eq:system_w_lumped}
	\bH_m c_m = - g_m, \; m = 1, \dots, N,
	\eeq
	where $g_m$ is given as in \eqref{eq:gradient_w}, and $\bH_m$ is given by the lumped Hessian
	\beq\label{eq:Hessian_m}
	\bH_m := \sum_{n=1}^N \bH_{mn}, \; m = 1, \dots, N,
	\eeq
	with $\bH_{mn}$ defined in \eqref{eq:Hessian_mn}.
	We refer to \cite{DetommasoCuiMarzoukEtAl18} for this technique and a diagonalization $\bH_m = \bH_{mm}$. Moreover, to find a good step size $\varepsilon$ in \eqref{eq:mapl_w}, we adopt a classical line search \cite{NocedalWright2006}, see Appendix B.
	
	\subsection{Hessian-based subspace}
	\label{sec:Hessian_subspace}
	To construct a data-informed subspace of the parameter space, we exploit the geometry of the posterior density characterized by its Hessian. More specifically, we seek the basis functions $\psi_i$, $i = 1, \dots, r$, as the eigenvectors corresponding to the $r$ largest eigenvalues of the generalized eigenvalue problem 
	\beq\label{eq:gEigenDecom}
	\bE[\nabla_x^2 \eta_y(x)] \psi_i = \lambda_i \Gamma_0^{-1} \psi_i, \quad i = 1, \dots, r, 
	\eeq
	where $\Gamma_0$ is the covariance of $x$ under the prior distribution (not necessarily Gaussian), $ \psi_i^T \Gamma_0^{-1} \psi_j = \delta_{ij}$, $i, j = 1, \dots, r$, $\bE[\nabla_x^2 \eta_y(x)]$, with $\eta_y(x) := -\log(p_\xi(y-f(x)))$, is the averaged Hessian of the negative log-likelihood function w.r.t.\ a certain distribution, e.g., the prior, posterior, or Gaussian approximate distribution \cite{CuiMarzoukWillcox2016}. 
	Here we propose to evaluate $\bE[\nabla_x^2 \eta_y(x)]$ by an adaptive sample average approximation at the samples pushed from the prior to the posterior, and adaptively construct the eigenvectors $\Psi$, as presented in next section.
	% Without loss of generality, we consider a Gaussian prior $\mu_0 = \cN(\bar{x}, C_{\text{pr}})$ with mean $\bar{x}$ and covariance $C_{\text{pr}} \in \bR^{d\times d}$, which is symmetric, positive, and definite. 
	For linear Bayesian inference problems, with $f(x) = Ax$ for $A \in \bR^{s \times d}$, a Gaussian prior distribution $x \sim \cN(\bar{x}, \Gamma_0)$ and a Gaussian noise $\xi \sim \cN(0, \Gamma_\xi)$ lead to a Gaussian posterior distribution given by $\cN(x_{\text{MAP}}, \Gamma_{\text{post}})$, where \cite{Stuart10}
	\beq\label{eq:Gauss_posterior}
	\Gamma^{-1}_{\text{post}} = \nabla_x^2\eta_y + \Gamma_0^{-1}, \; x_{\text{MAP}} = \bar{x} - \Gamma_{\text{post}} A^T\Gamma_\xi^{-1} (y-A\bar{x}).
	\eeq
	% It is well known that the eigenvalues $\lambda_i$, $i = 1, 2, \dots, $ of the generalized Hermetian eigenvalue problem 
	% \beq
	% \nabla^2_x \eta_y \psi_i = \lambda_i C_{\text{pr}}^{-1} \psi_i, \text{ where } \psi_i^{\top} C_{\text{pr}}^{-1} \psi_j  = \triangle_{ij},
	% \eeq
	Therefore, the eigenvalue $\lambda_i$ of $(\nabla_x^2\eta_y, \Gamma_0^{-1})$, with $\nabla_x^2\eta_y = A^T \Gamma_\xi^{-1} A$, measures the relative variation between the data-dependent log-likelihood and the prior in direction $\psi_i$. For $\lambda_i \ll 1$, 
	% the data result in negligible variation of the posterior compared to the prior in direction $\psi_i$, or in another words, 
	the data provides negligible information in direction $\psi_i$, so the difference between the posterior and the prior in $\psi_i$ is negligible. In fact, it is shown in \cite{SpantiniSolonenCuiEtAl15} that the subspace constructed by \eqref{eq:gEigenDecom} is optimal for linear $f$.  
	% For general nonlinear Bayesian inference problems where the posterior is not necessarily Gaussian, the data-informed directions can be similarly obtained by the eigenvectors corresponding to the largest eigenvalues of the generalized Hermetian eigenvalue problem, 
	% \beq\label{eq:gEigenDecom}
	% H \psi_i = \lambda_i C_{\text{pr}}^{-1} \psi_i, \text{ where } \psi_i^{\top} C_{\text{pr}}^{-1} \psi_j  = \triangle_{ij},
	% \eeq
	% where $H$ can be taken as an averaged Hessian \cite{CuiLawMarzouk16}
	% \beq\label{eq:Hessian_average}
	% H = \bE^{\mu_y}[\nabla^2_x \eta_y] \approx \frac{1}{N} \sum_{n = 1}^N \nabla^2_x \eta_y(x_n),
	% \eeq
	% or a combined Hessian \cite{chenGhattas18b} to account for the variation of the Hessian at different samples $x_n$, $n = 1, \dots, N$, from the posterior distribution.
	Let $(\lambda_i, \psi_i)_{1 \leq i \leq r}$ denote the $r$ largest eigenpairs such that $|\lambda_1| \geq |\lambda_2| \geq \cdots \geq |\lambda_r| \geq \varepsilon_\lambda > |\lambda_{r+1}|$ for some small tolerance $\varepsilon_\lambda < 1$. Then the Hessian-based subspace spanned by the eigenvectors $\Psi = (\psi_1, \dots, \psi_r)$ captures the most variation of the parameter $x$ informed by data $y$.
	We remark that to solve the generalized Hermitian eigenvalue problem \eqref{eq:gEigenDecom}, we employ a randomized SVD algorithm \cite{HalkoMartinssonTropp11}, which requires $O(NrC_h + d r^2)$ flops, where $C_h$ is the cost of a Hessian action in a direction.

	\subsection{Parallel and adaptive pSVN algorithm}
	\label{sec:adaptivepSVN}
	Given the bases $\Psi$ as the data-informed parameter directions, we can draw samples $x_1, \dots, x_N$ from the prior distribution and push them by pSVN to match the posterior distribution in a low-dimensional subspace, while keeping the components of the samples in the complementary subspace unchanged. We set the stopping criterion as: (i) the maximum norm of the updates $w_m^{l} - w_m^{l-1}$, $m = 1, \dots, N$, is smaller than a given tolerance $\text{Tol}_g$; (ii) the maximum norm of the gradients $g_m$, $m = 1, \dots, N$, is smaller than a given tolerance $\text{Tol}_w$; or (iii) the number of iterations $l$ reaches a preset number $L$.
	Moreover, we take advantage of pSVN advantages in low-dimensional subspaces---including fast computation, lightweight communication, and low memory footprint---and provide an efficient parallel implementation using MPI communication
	%of the adaptive pSVN algorithm. 
	%and compare it to a parallel SVN algorithm without using the projection.
	% We present a parallel pSVN using MPI communication 
	in Algorithm \ref{alg:pSVN_parallel}, with analysis in Appendix C.

	% \begin{algorithm}[!htb]
	%   \caption{pSVN}
	%   \label{alg:pSVN}
	% \begin{algorithmic}[1]
	%   \STATE {\bfseries Input:} prior samples $x_1, \dots, x_N$, bases $\Psi$, density $p_y$.
	%   \STATE {\bfseries Output:} posterior samples $x^y_1, \dots, x_N^y$.
	%   \STATE Perform projection \eqref{eq:projection} with $x_n = x_n^r + x_n^\perp$ and the samples $w^{l-1}_n$, $n = 1, \dots, N$, at $l=1$.
	%   \REPEAT 
	%       \STATE Compute the gradient and Hessian by \eqref{eq:gradient_w} and \eqref{eq:Hessian_w}.
	%       \STATE Compute the kernel and its gradient by \eqref{eq:kernel} and \eqref{eq:metric}.
	%       \STATE Assemble and solve system \eqref{eq:system_w_lumped} for $c_1, \dots, c_N$.
	%       \STATE Perform a line search for $\varepsilon$ to get $w_1^l, \dots, w_N^l$ by \eqref{eq:mapl_w}.
	%       \STATE Update the samples $x^r_n = \Psi w^l_n + \bar{x}$, $n =1, \dots, N$.
	%       \STATE Set $l \leftarrow l + 1$.
	%   \UNTIL{A stopping criterion is met.}
	%   \STATE Reconstruct samples $x_n^y = x_n^r + x^\perp$, $n = 1, \dots, N$.
	% \end{algorithmic}
	% \end{algorithm}

	% This is presented in Algorithm \ref{alg:pSVN_two-level}. We remark that the same stopping criteria in Algorithm \ref{alg:pSVN_two-level} as those in Algorithm \ref{alg:pSVN} are used with smaller tolerances $\text{Tol}_g^2, \text{Tol}_w^2$ for the gradients and the updates, e.g., $\text{Tol}_g^2 = 10^{-1}\text{Tol}_g$, and $\text{Tol}_w^2 = 10^{-1}\text{Tol}_w$.

	% \subsection{Parallel computation and implementation}
	
	\begin{algorithm}[!htb]
		\caption{pSVN in parallel using MPI}
		\label{alg:pSVN_parallel}
		\begin{algorithmic}[1]
			\STATE {\bfseries Input:} $M$ prior samples, $x_1, \dots, x_M$, in each of $K$ cores, bases $\Psi$, and density $p_y$ in all cores.
			\STATE {\bfseries Output:} posterior samples $x^y_1, \dots, x_M^y$ in each core.
			\STATE Perform projection \eqref{eq:projection} to get $x_m = x_m^r + x_m^\perp$ and the samples $w^{l-1}_m$, $m = 1, \dots, M$, at $l=1$.
			\STATE \emph{\color{black} Perform \emph{MPI\_Allgather} for $w^{l-1}_m$, $m = 1, \dots, M$.}
			\REPEAT 
			\STATE Compute the gradient and Hessian by \eqref{eq:gradient_w}.
			\STATE \emph{\color{black}  Perform \emph{MPI\_Allgather} for the gradient and Hessian.}
			\STATE Compute the kernel and its gradient by \eqref{eq:kernel} and \eqref{eq:metric}.
			\STATE \emph{\color{black}  Perform \emph{MPI\_Allgather} for $k_m$, $m = 1, \dots, M$, \\ \emph{MPI\_Allreduce} w. sum for $\sum_m k_m$ and $\sum_m \nabla_w k_m$.}
			\STATE Assemble and solve system \eqref{eq:system_w_lumped} for $c_1, \dots, c_M$.
			\STATE Perform a line search to get $w_1^l, \dots, w_M^l$.
			\STATE \emph{\color{black} Perform \emph{MPI\_Allgather} for $w^{l}_m$, $m = 1, \dots, M$.}
			\STATE Update the samples $x^r_m = \Psi w^l_m + \bar{x}$, $m =1, \dots, M$.
			\STATE Set $l \leftarrow l + 1$.
			\UNTIL{A stopping criterion is met.}
			\STATE Reconstruct samples $x_m^y = x_m^r + x^\perp_m$, $m = 1, \dots, M$.
		\end{algorithmic}
	\end{algorithm}
	
	In Algorithm \ref{alg:pSVN_parallel}, we assume that the bases $\Psi$ for the projection are the data informed parameter directions, which are obtained by the Hessian-based algorithm in Section \ref{sec:Hessian_subspace} at the ``representative'' samples $x_1, \dots, x_N$. However, we do not have these samples but only the prior samples at the beginning. To address this problem, we propose an adaptive algorithm that adaptively construct the bases $\Psi$ based on samples pushed forward from the prior to the posterior, see Algorithm \ref{alg:pSVN_two-level}.
	
	\begin{algorithm}[!htb]
		\caption{Adaptive pSVN}
		\label{alg:pSVN_two-level}
		\begin{algorithmic}[1]
			\STATE {\bfseries Input:} $M$ prior samples, $x_1, \dots, x_M$, in each of $K$ cores, and density $p_y$ in all cores.
			\STATE {\bfseries Output:} posterior samples $x^y_1, \dots, x_M^y$ in each core.
			\STATE Set level $l_2 = 1$, $x_m^{l_2-1} = x_m$, $m = 1, \dots, M$.
			\REPEAT 
			\STATE Perform the eigendecomposition \eqref{eq:gEigenDecom} at samples $x_1^{l_2-1}, \dots, x_M^{l_2-1}$, and form the bases $\Psi^{l_2}$.
			\STATE Apply \textbf{Algorithm} \ref{alg:pSVN_parallel} to update the samples \\ $[x_1^{l_2}, \dots, x_M^{l_2}] = \text{pSVN}([x_1^{l_2-1}, \dots, x_M^{l_2-1}], K, \Psi^{l_2}, p_y)$.
			\STATE Set $l_2 \leftarrow l_2 + 1$.
			\UNTIL{A stopping criterion is met.}
		\end{algorithmic}
	\end{algorithm}
	
	\section{Numerical experiments}
	\label{sec:numerics}
	We demonstrate the convergence, accuracy, and dimension-independence of the pSVN method by two examples, one a linear problem with Gaussian posterior to demonstrate the convergence and accuracy of pSVN in comparison with SVN and SVGD, the other a nonlinear problem to demonstrate accuracy as well as the dimension-independent and sample-independent convergence of pSVN and its scalability w.r.t.\ the number of processor cores. The code is described in Appendix D.
	%For the latter purpose, we also consider a Bayesian autoencoder problem. 
	% The code is available at our Bitbucket repository \cite{pSVN} for the test examples. \pc{add link}
	%and \cite{pSVNnn} for the third.
	
	\subsection{A linear inference problem}
	For the linear inference problem, we have the parameter-to-observable map
	\beq
	f(x) = A x,  
	\eeq
	where the linear map $A = O (B x)$, with an observation map $O:\bR^d \to \bR^s$, and an inverse discrete differential operator $B = (L + M)^{-1} :\bR^d \to \bR^d$ where $L$ and $M$ are the discrete Laplacian and mass matrices in the PDE model $ -\triangle u + u = x,  \text{ in } (0, 1), \; u(0)=0, \; u(1)=1.$
	$s = 15$ pointwise observations of $u$ with $1\%$ noise are distributed with equal distance in $(0, 1)$. The input $x$ is a random field with Gaussian prior $\cN(0, \Gamma_0)$, where $\Gamma_0$ is discretized from $(I - 0.1 \triangle)^{-1}$ with identity $I$ and Laplace operator $\triangle$. We discretize this forward model by a finite element method with piecewise linear elements on a uniform mesh of size $2^n$, which leads to the parameter dimension $d = 2^n + 1$. 
	
	\begin{figure}[!htb]
		\begin{center}
			\includegraphics[width=0.32\columnwidth]{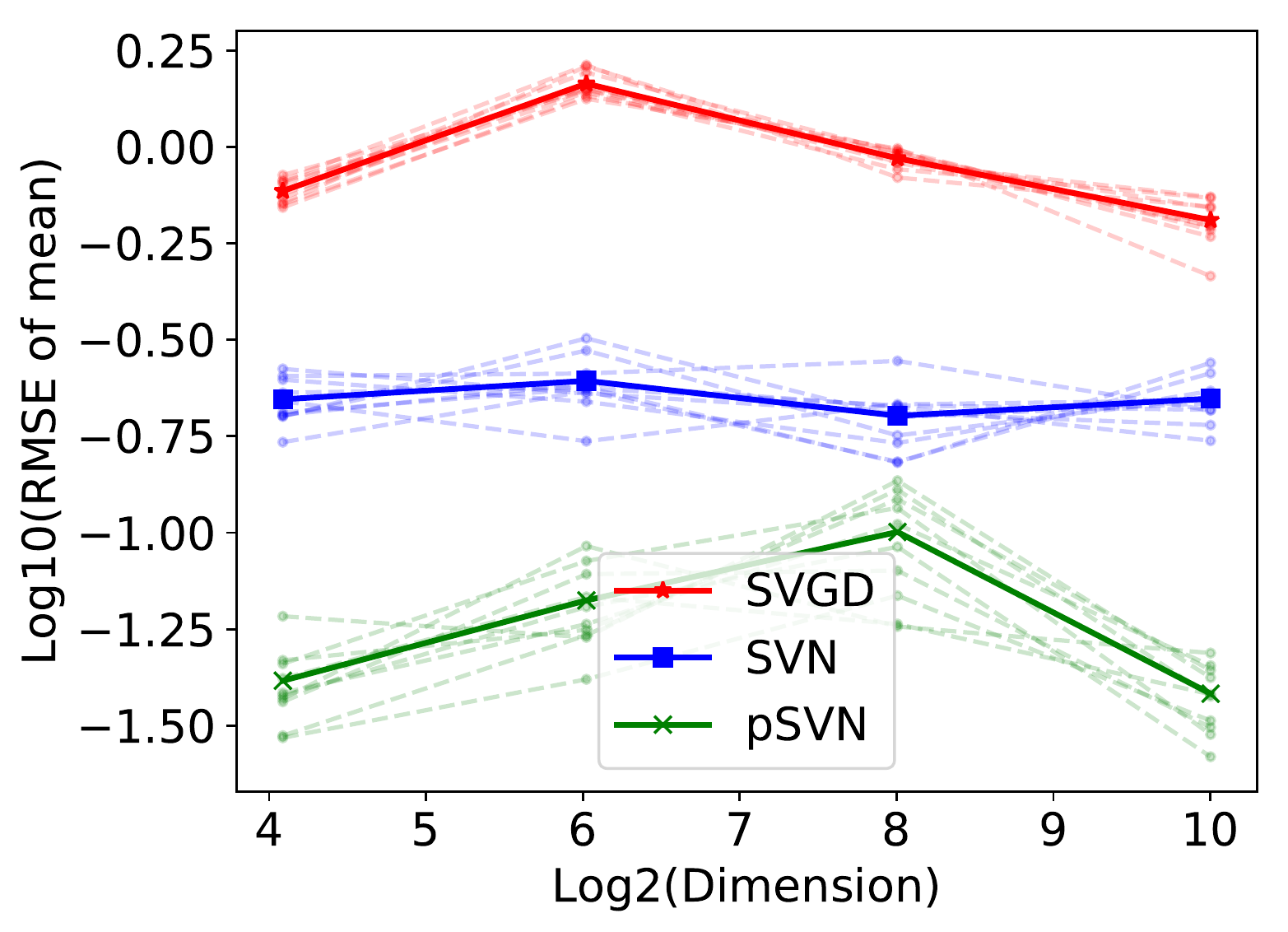}
			\includegraphics[width=0.32\columnwidth]{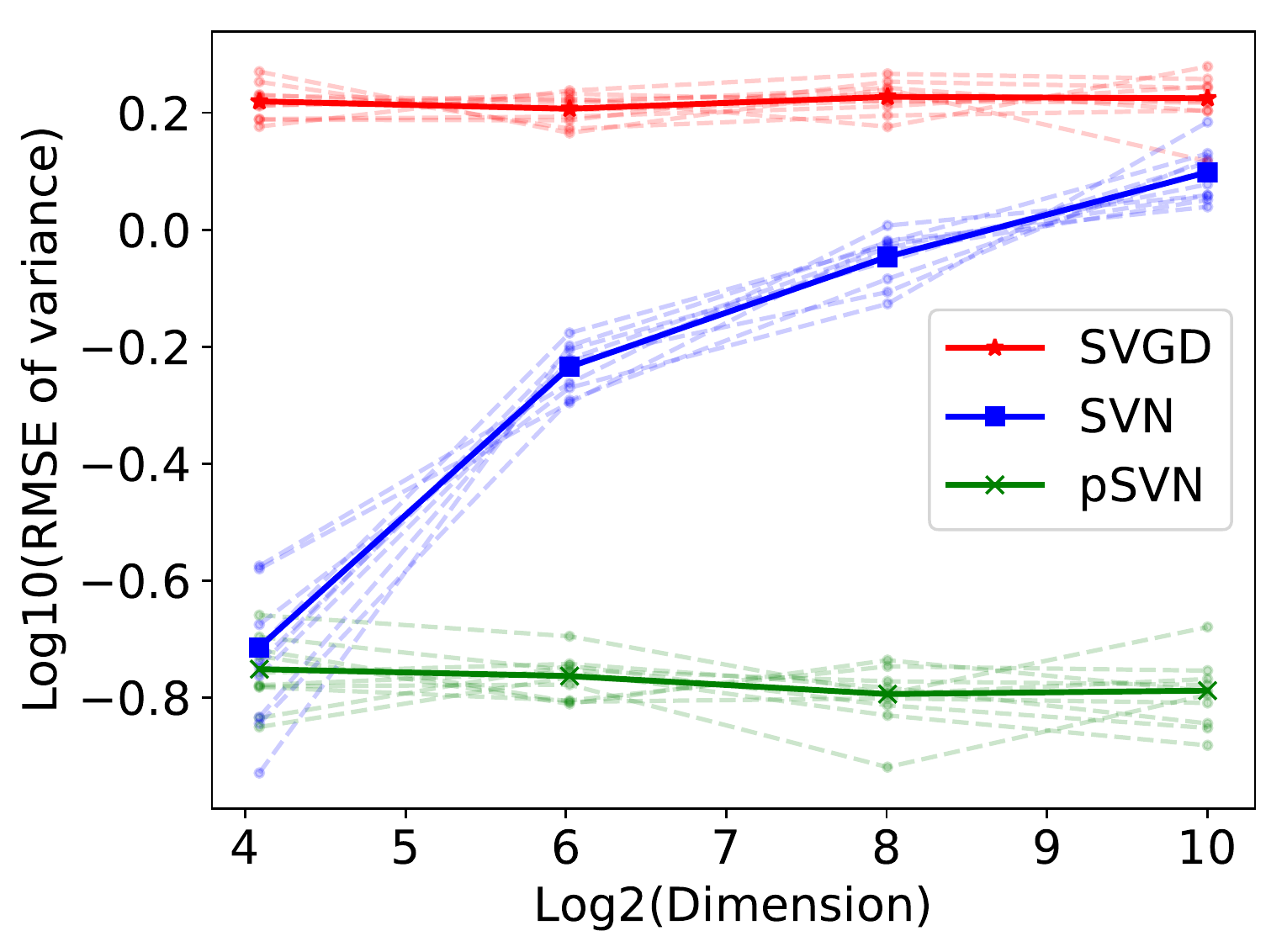}
			\includegraphics[width=0.32\columnwidth]{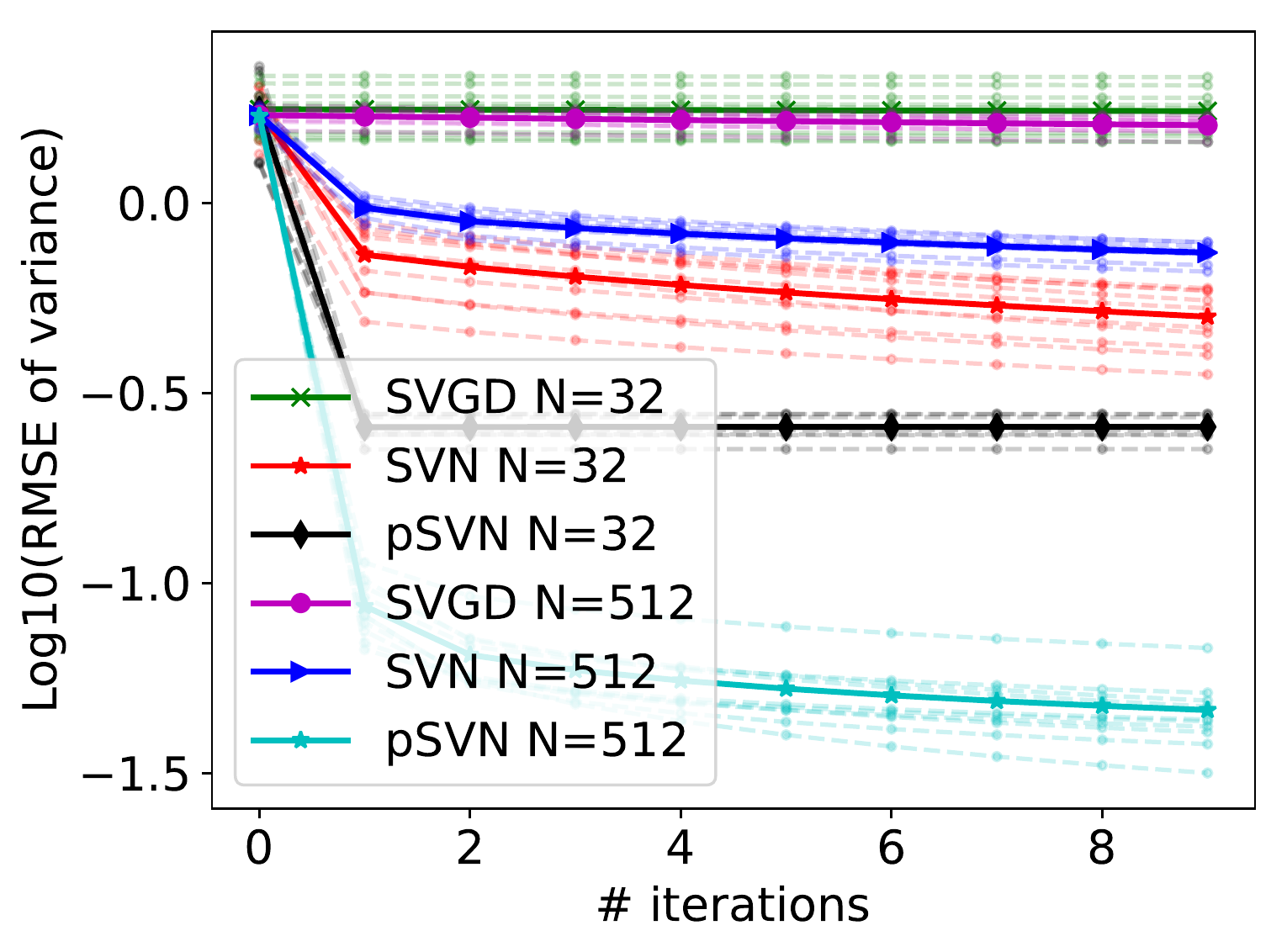}
			\vspace{-0.4cm}
			\caption{Decay of the RMSE (with 10 trials in dashed lines) of the L2-norm of the mean (left) and pointwise variance (middle) of the parameter w.r.t.\ dimension $d = 16, 64, 256, 1024$ with $N = 128$ samples. Right: Decay of the RMSE of the L2-norm of the pointwise variance with $N = 32, 512$ samples in parameter dimension $d = 256$ w.r.t.\ \# iterations. Comparison for SVGD, SVN, pSVN.}
			\label{Dimension}
		\end{center}
		\vspace{-0.4cm}
	\end{figure}

	% \begin{figure}[!htb]
	% \begin{center}
	% \includegraphics[width=0.45\columnwidth]{figure/error_mean_MSE}
	% % \vspace{-0.4cm}
	% % \caption{Decay of the mean squared error of the mean of the parameter using different number of samples for pSVN vs SVN.}
	% % \label{MSE-mean}
	% % \end{center}
	% % \end{figure}
	% % \begin{figure}[!htb]
	% % \begin{center}
	% \includegraphics[width=0.45\columnwidth]{figure/error_variance_MSE}
	% \vspace{-0.4cm}
	% \caption{Decay of the mean squared errors (with 10 trials) of the L2-norm of the mean (top) and pointwise variance (bottom) of the parameter using 32 and 512 samples by SVGD, SVN, and pSVN.}
	% \label{MSE}
	% \end{center}
	% \vspace{-0.4cm}
	% \end{figure}
	
	% The L2-norm of the mean and pointwise variance of the parameter $x$ w.r.t.\ its posterior distribution can be explicitly computed using \eqref{eq:Gauss_posterior}, which serve as the reference for the sample approximation errors. 
	Figure \ref{Dimension} compares the convergence and accuracy of SVGD, SVN, and pSVN by the decay of the root mean square errors (RMSE) (using 10 trials and 10 iterations) of the sample mean and variance (with L2-norm of errors computed against analytic values in \eqref{eq:Gauss_posterior}) w.r.t.\ parameter dimensions and iterations. We  observe much faster convergence and greater accuracy of pSVN relative to SVGD and SVN, for both mean and especially variance, which measures the goodness of samples. In particular, we see from the middle figure that the SVN estimate of variance deteriorates quickly with increasing dimension, while pSVN leads to equally good variance estimate. Moreover, from the right figure we can see that pSVN converges very rapidly in a subspace of dimension 6 (at tolerance $\varepsilon_\lambda = 0.01$ in Section \ref{sec:Hessian_subspace}, i.e., $|\lambda_7| < 0.01$) and achieves higher accuracy with larger number of samples, while SVN converges slowly and leads to large errors. With the same number of iterations of SVN and pSVN, SVGD produces no evident error decay.  
	
	\subsection{A nonlinear inference problem}
	We consider a nonlinear benchmark inference problem (which is often used for testing high-dimensional inference methods \cite{Stuart10, CuiLawMarzouk16, BESKOS2017327}), whose forward map is given by 
	$f(x) = O(S(x))$, with observation map $O:\bR^d \to \bR^s$ and a nonlinear solution map $u = S(x) \in \bR^d$ of the lognormal diffusion model $ -\nabla\cdot (e^x \nabla u) = 0,  \text{ in } (0, 1)^2$ with $u = 1$ on top and $u=0$ on bottom boundaries, and zero Neumann conditions on left and right boundaries. 
	$49$ pointwise observations of $u$ are equally distributed in $(0, 1)^2$. We use $10\%$ noise to test accuracy against a DILI MCMC method \cite{CuiLawMarzouk16} with 10,000 MCMC samples as reference and a challenging $1\%$ noise for a dimension-independence test of pSVN. The input $x$ is a random field with Gaussian prior $\cN(0, \Gamma_0)$, where $\Gamma_0$ is a discretization of $(I - 0.1 \triangle)^{-2}$. We solve this forward model by a finite element method with piecewise linear elements on a uniform mesh of varying sizes, which leads to a sequence of parameter dimensions. 
	
	\begin{figure}[!htb]
		\begin{center}
			\includegraphics[width=0.35\columnwidth]{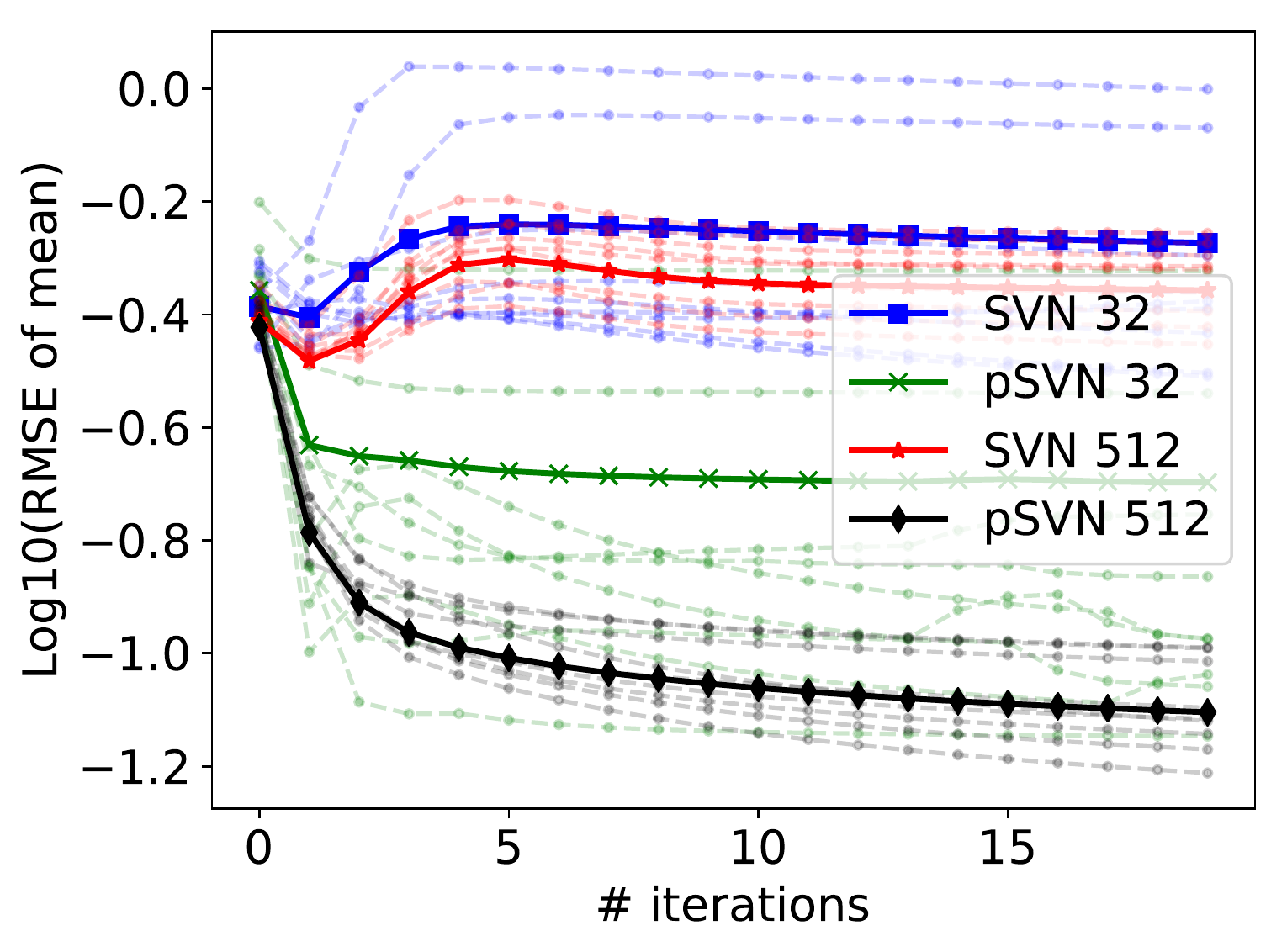}
			\includegraphics[width=0.35\columnwidth]{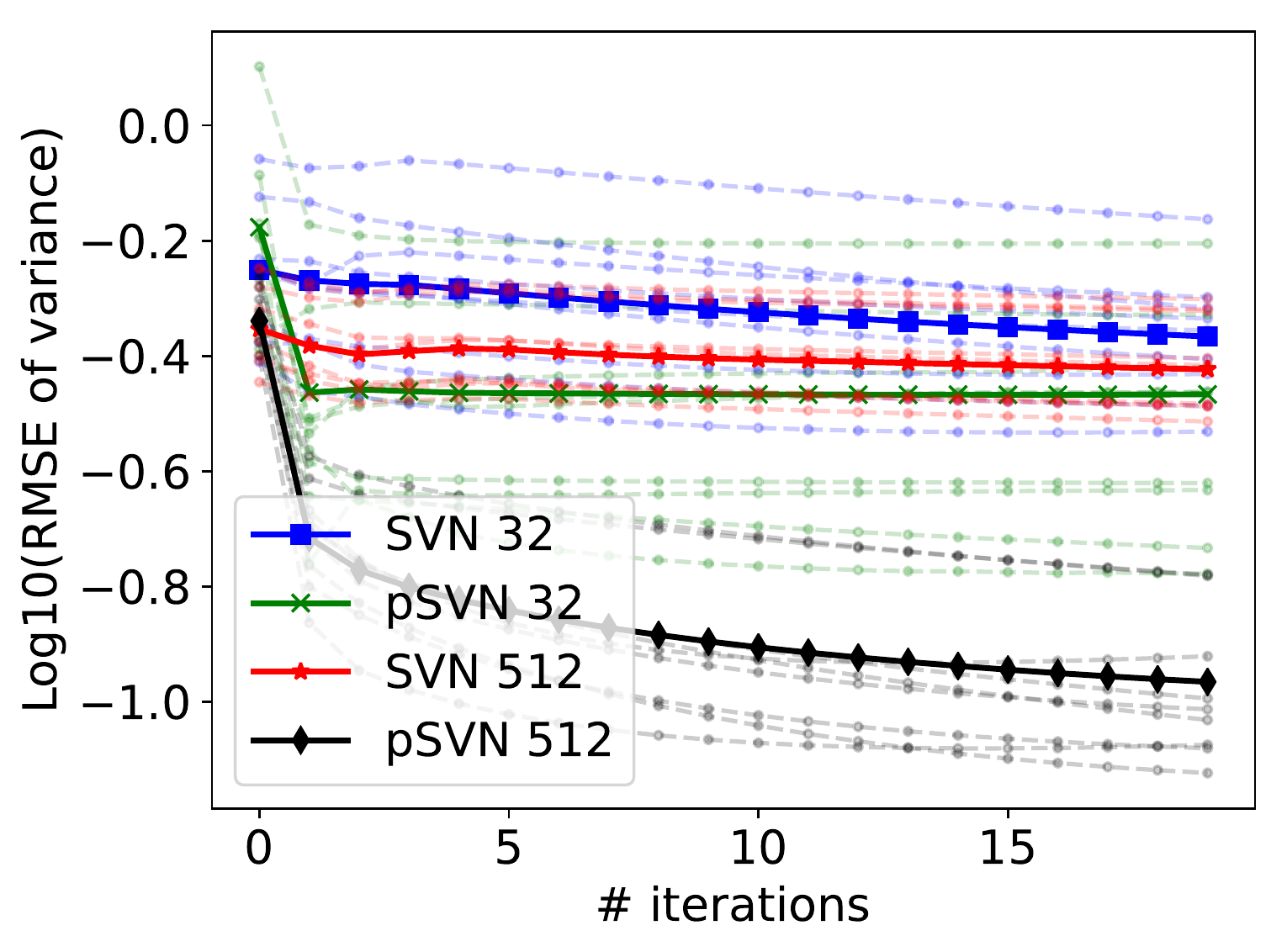}
			% \includegraphics[width=0.32\columnwidth]{figure/core_time}
			% \vspace{-0.4cm}
			\caption{Decay of the RMSE (with 10 trials in dashed lines) of the L2-norm of the mean (left) and pointwise variance (right) of the parameter with dimension $d = 1089$ and $N = 32, 512$ samples.}
			\label{accuracy}
		\end{center}
	\end{figure}
	
	Figure \ref{accuracy} shows the comparison of the accuracy and convergence of pSVN and SVN for their sample estimate of mean and variance. We can see that in high dimension, $d = 1089$, pSVN converges faster and achieves higher accuracy than SVN for both mean and variance estimate. Moreover, SVN using the kernel \eqref{eq:kernel} in high dimensions (involving low-rank decomposition of the metric $M$ for high-dimensional nonlinear problems) is more expensive than pSVN per iteration.

	\begin{figure}[!htb]
		\begin{center}
			\includegraphics[width=0.32\columnwidth]{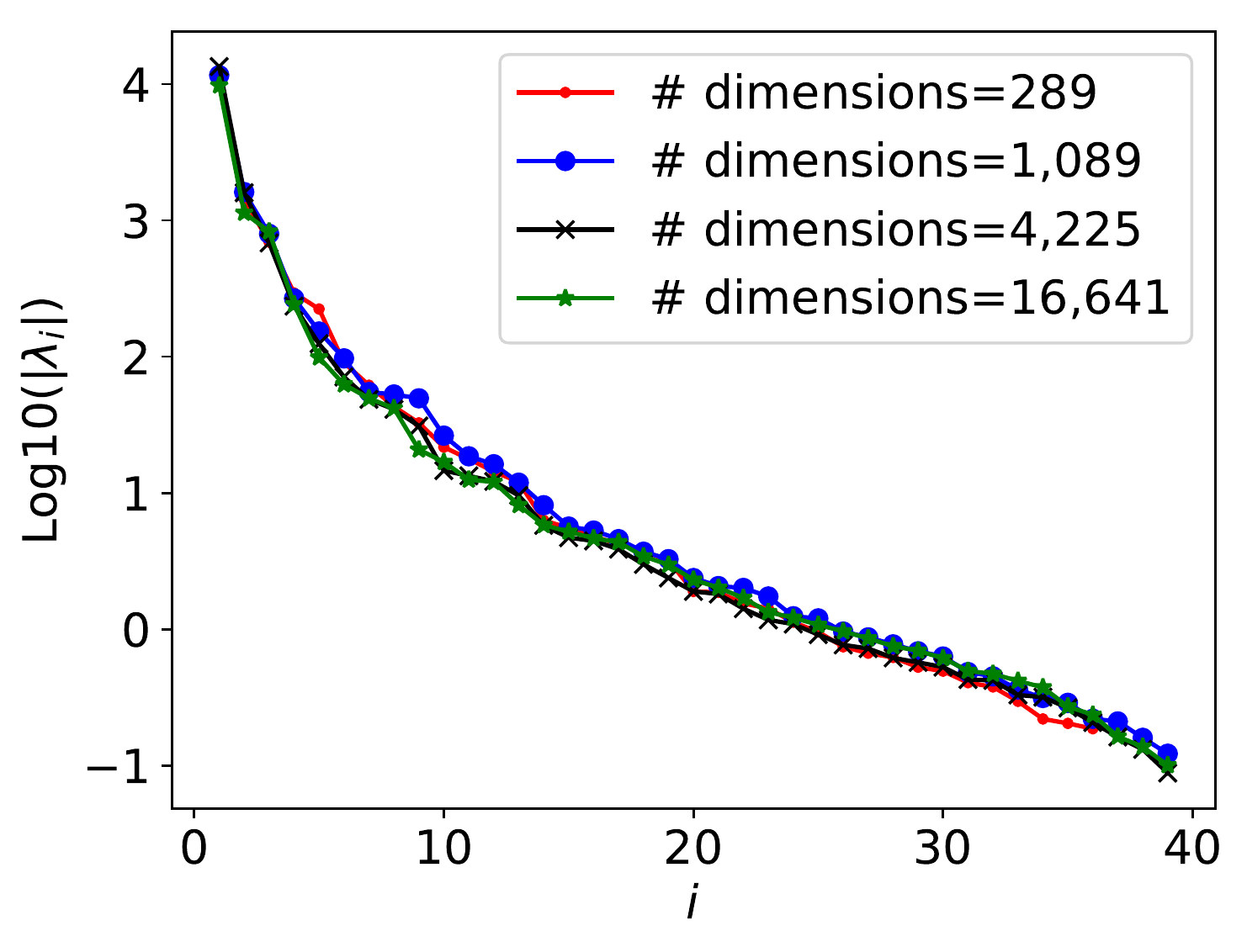}
			\includegraphics[width=0.33\columnwidth]{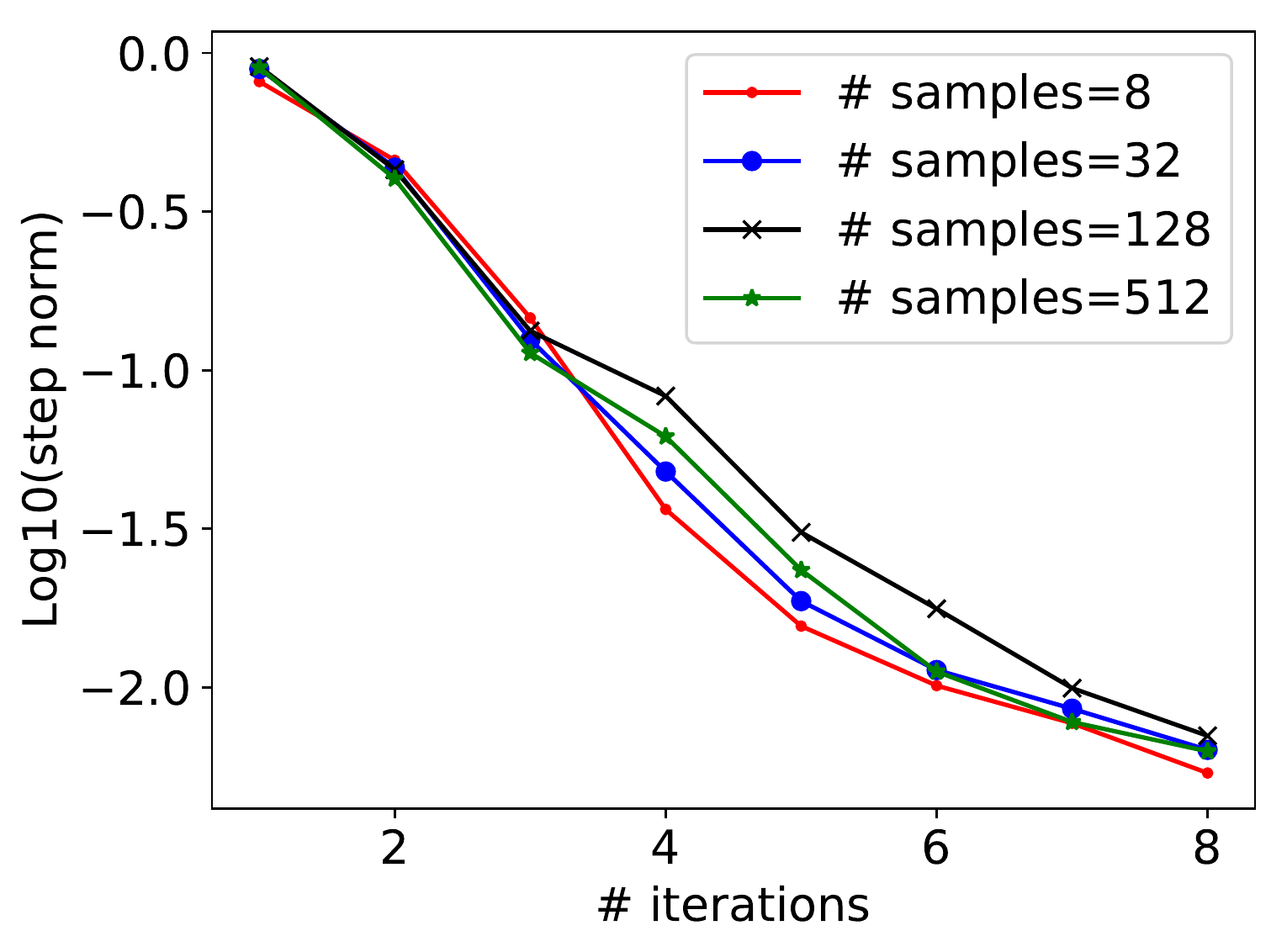}
			\includegraphics[width=0.32\columnwidth]{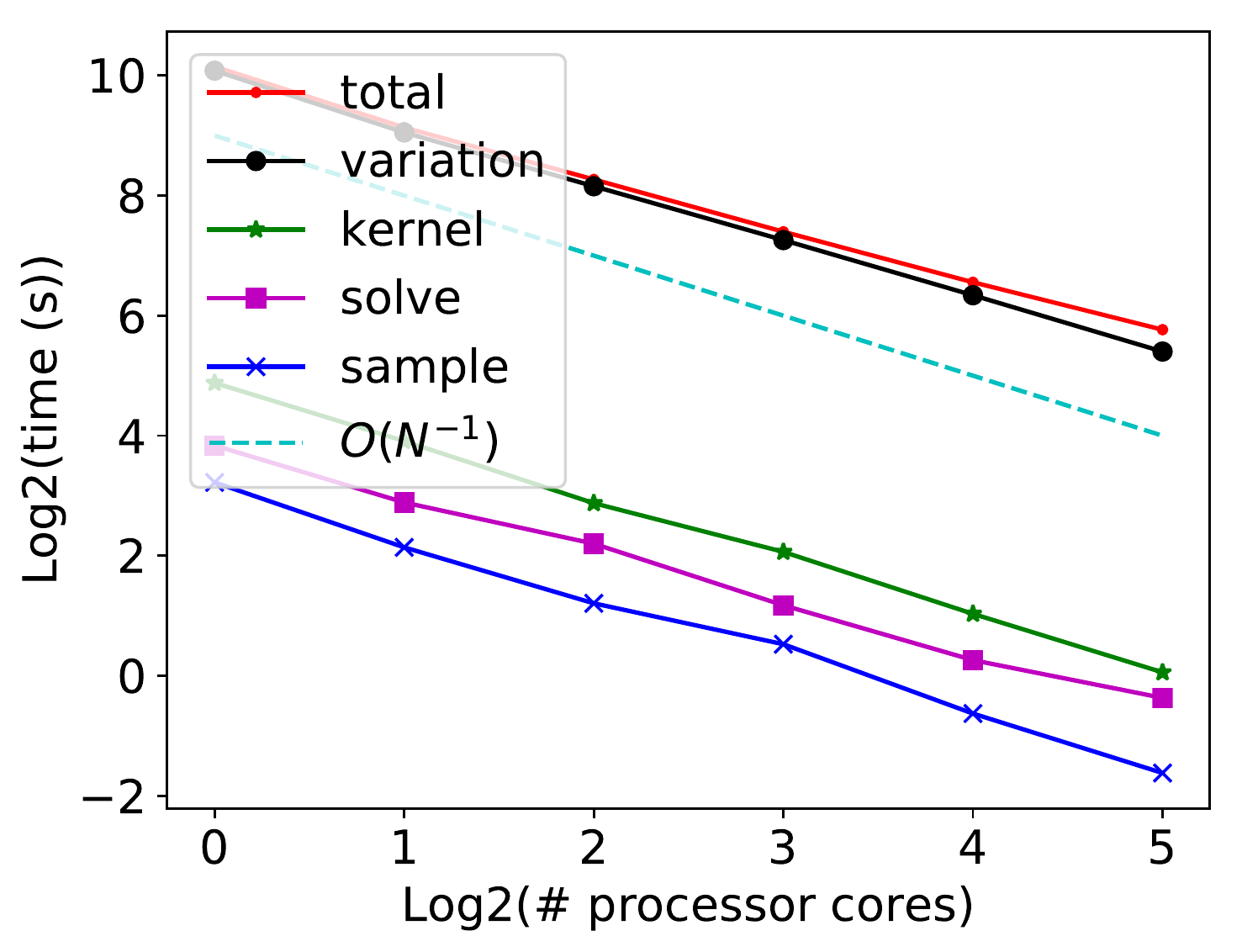}
			\vspace{-0.4cm}
			\caption{Left: Decay of eigenvalues $\log_{10}(|\lambda_i|)$ with increasing dimension $d$. Middle: Decay of a stopping criterion---the averaged norm of the update $w^l- w^{l-1}$ w.r.t.\ the iteration number $l$, with increasing number of samples. Right: Decay of the wall clock time (seconds) of different computational components w.r.t.\ increasing number of processor cores on a log-log scale.}
			\label{scalability-dimensions}
		\end{center}
	\end{figure}
	
	We next demonstrate pSVN's independence of the number of parameter and sample dimensions, and its scalability w.r.t.\ processor cores. First, the dimension of the Hessian-based subspace $r$, which determines the computational cost of pSVN, depends on the decay of the absolute eigenvalues $|\lambda_i|$ as presented in Section \ref{sec:Hessian_subspace}. The left part of Figure \ref{scalability-dimensions} shows that with increasing $d$ from 289 to over 16K, $r$ does not change, which implies that the convergence of pSVN is independent of the number of nominal parameter dimensions. Second, as shown in the middle part of Figure \ref{scalability-dimensions}, with increasing number of samples for a fixed parameter dimension $d = 1089$, the averaged norm of the update $w^l- w^{l-1}$, as one convergence indicator presented in Subsection \ref{sec:adaptivepSVN}, decays similarly, which demonstrates the independence of the convergence of pSVN w.r.t.\ the number of samples.
	% \begin{figure}[!htb]
	% \begin{center}
	% \includegraphics[width=0.8\columnwidth]{figure/step_norm}
	% \vspace{-0.4cm}
	% \caption{Decay of the averaged norm of the update $w^l- w^{l-1}$ w.r.t.\ the iteration number $l$, with increasing number of samples.}
	% \label{scalability-samples}
	% \end{center}
	% \end{figure}
	Third, in the right of Figure \ref{scalability-dimensions} we plot the \textit{total} wall clock time of pSVN and the time for its computational components in Algorithm \ref{alg:pSVN_parallel} using different number of processor cores for the same work, i.e., the same number of samples (256), including \textit{variation} for forward model solve, gradient and Hessian evaluation, as well as eigendecomposition, \textit{kernel} for kernel and its gradient evaluation, \textit{solve} for solving the Newton system \eqref{eq:system_w_lumped}, and \textit{sample} for sample projection and reconstruction. We can observe nearly perfect strong scaling w.r.t.\ increasing number of processor cores. Moreover, the time for \textit{variation}, which depends on parameter dimension $d$, dominates the time for all other components, in particular \textit{kernel} and \textit{solve} whose cost only depends on $r$, not $d$.

	\section{Conclusion}
	We presented a fast and scalable variational method, 
	pSVN, for Bayesian inference in high dimensions. The method exploits the geometric structure of the posterior via its Hessian, and the intrinsic low-dimensionality of the change from prior to posterior characteristic of many high-dimensional inference problems via low rank approximation of the averaged Hessian of the log likelihood, computed efficiently using randomized matrix-free SVD. The fast convergence and higher accuracy of pSVN relative to SVGD and SVN, its complexity that is independent of parameter and sample dimensions, and its scalability w.r.t.\ processor cores were demonstrated for linear and nonlinear inference problems. Investigation of pSVN to tackle intrinsically high-dimensional inference problem (e.g., performed in local dimensions as the message passing scheme or combined with dimension-independent MCMC to update samples in complement subspace) is ongoing. Further development and application of pSVN to more general probability distributions, projection basis constructions, and forward models such as deep neural network, and further analysis of the convergence and scalability of pSVN w.r.t.\ the number of samples, parameter dimension reduction, and data volume, are of great interest.

	\bibliographystyle{plain}
	\bibliography{pSVN}

\newpage

\section*{Appendix A: Proof of Theorem \ref{thm:convergence}}

\begin{proof}
	By definition of the posterior density $p(x|y)$ in \eqref{eq:Bayes_x} and the projected posterior density $p^r(x|y)$ in \eqref{eq:Bayes_w}, we have 
	\beq\label{eq:KLsplit}
	\begin{split}
		\cD_{\text{KL}}(p(x|y) \,|\, p^r(x|y)) & = \int_{\bR^d}  \log\left(\frac{p_y(x)}{p_y^r(x)} \frac{Z_r}{Z}\right) \frac{1}{Z} p_y(x) dx \\
		& = \int_{\bR^d} (\eta_y(x^r) - \eta_y(x))  \frac{1}{Z} p_y(x) dx  + \log\left(\frac{Z_r}{Z}\right), 
	\end{split}
	\eeq
	where we used the definitions of $p_y(x)$ and $p^r_y(x)$ in \eqref{eq:Bayes_x} and \eqref{eq:Bayes_w} in the second equality.
	By definition of $\eta_y$ in \eqref{eq:potential}, we have \beq\label{eq:boundeta}
	\begin{split}
		\eta_y(x^r) - \eta_y(x) & = \frac{1}{2} ||y - f(x^r)||_\Gamma^2 - \frac{1}{2} ||y - f(x)||_\Gamma^2\\
		& = y^T \Gamma^{-1} (f(x) - f(x^r)) - \frac{1}{2} (f(x)+f(x^r))^T \Gamma^{-1} (f(x) - f(x^r)) \\
		& \leq ||y^T||_\Gamma ||f(x) - f(x^r)||_\Gamma + \frac{1}{2} ||f(x) + f(x^r)||_\Gamma ||f(x) - f(x^r)||_\Gamma\\
		& \leq \frac{C_b}{2}  (2||y^T||_\Gamma + ||f(x)||_\Gamma + ||f(x^r)||_\Gamma )||x - x^r||_X 
	\end{split}
	\eeq
	where we used Assumption \ref{ass:f} in the second inequality for $\max\{||x||_X, ||x^r||_X\} < b$. 
	Therefore, the first integral in \eqref{eq:KLsplit}, denoted as (I) can be bounded by (note that $\exp(-\eta_y(\cdot)) \leq 1$)
	\beq
	\begin{split}
		(I) \leq \frac{C_b}{2Z} \int_{\bR^d} \left(2||y^T||_\Gamma  + ||f(x)||_\Gamma + ||f(x^r)||_\Gamma  \right) p_0(x) dx \; ||x - x^r||_X,
	\end{split}
	\eeq
	
	By Assumption \ref{ass:f}, we have 
	\beq
	(I) \leq C_I ||x - x^r||_X,
	\eeq
	for a constant $C_I = C_b (||y^T||_\Gamma + C_f)/Z$. 
	
	For the second term $\log(Z_r/Z)$ in \eqref{eq:KLsplit}, we have for 
	\beq
	\begin{split}
		\Big| 1 - \frac{Z_r}{Z} \Big| & = \frac{1}{Z} |Z - Z_r| \\
		&  \leq  \frac{1}{Z} \int_{\bR^d} |\exp(-\eta_y) - \exp(-\eta_y^r)  | p_0(x) dx\\
		& \leq \frac{1}{Z} \int_{\bR^d}  |\eta_y - \eta_y^r|  p_0(x) dx \\
		& \leq C_I  ||x - x^r||_X,
	\end{split}
	\eeq
	where in the second inequality we used that $|e^{-\tau_1}-e^{-\tau_2}| < |\tau_1 - \tau_2|$ for $\tau_1, \tau_2 > 0$, for the last inequality we used the bound of the first integral of \eqref{eq:KLsplit}. Then by $\log(1+\tau) \leq \tau $ for $\tau \geq 0$, we have 
	\beq
	\log\left(\frac{Z_r}{Z}\right) \leq \log\left(1 + \Big| \frac{Z_r}{Z} - 1\Big|  \right) \leq  \Big| \frac{Z_r}{Z} - 1\Big| \leq C_I  ||x - x^r||_X,
	\eeq
	which completes the proof with constant $C = 2 C_I$
	
\end{proof}

\section*{Appendix B: Globalization by line search}
\label{sec:line-search}
%The step size $\varepsilon$ in \eqref{eq:mapl_w} is set as $1$ in the SVN algorithm in \pc{cite}, and as a small value, e.g., $10^{-2}$, in the SVGD \pc{cite}. However, 
Except for in the case of a linear inference problem, the cost functional---Kullback--Leibler divergence---is nonconvex. In the case of that the Newton approximation to the Kullback--Leibler divergence is locally exact, the simple choice of $\varepsilon=1$ is the optimal choice for the step size. However, since the geometry generally exhibits complex non-quadratic local structure, a constant stepsize $\varepsilon$ renders minimization of $\mathcal{D}_{\text{KL}}$ inefficient. A careful choice of the step size $\varepsilon$ is crucial for both fast convergence and stability of Stein variational methods. While, there are many options to choose from, we employ an Armijo line search globalization method to choose this step size, to much success. Specifically, at step $l =1, 2, \dots$, we seek $\varepsilon$ to minimize the Kullback--Leibler divergence 
\beq\label{eq:DKL_w}
\mathcal{D}_{\text{KL}} ((T_l)_* \pi_{l-1} | \pi_y) = \mathcal{D}_{\text{KL}} (\pi_{l-1} | (T_l)^*\pi_y),
\eeq
where $(T_l)^*$ is the pullback operator. Because 
\beq
\begin{split}
	& \mathcal{D}_{\text{KL}} (\pi_{l-1} | (T_l)^*\pi_y) = \bE_{\pi_{l-1}}[\log(\pi_{l-1}(\cdot))]\\
	& - \bE_{\pi_{l-1}}[\log(\pi_y(T_l(\cdot)) |\text{det} \; \nabla_w T_l(\cdot)|)],
\end{split}
\eeq
where the first term does not depend on $\varepsilon$. Hence we only need to consider the second term denoted as $\mathcal{D}_{\text{KL}}^{(2)}$, which is evaluated by the sample average approximation as
\beq\label{eq:DKL2}
\begin{split}
	\mathcal{D}_{\text{KL}}^{(2)} \approx & - \frac{1}{N} \sum_{n=1}^N \log(\pi_y(T_l(w_n^{l-1}))) \\
	& - \frac{1}{N} \sum_{n=1}^N \log(|\text{det} \; \nabla_w T_l(w_n^{l-1})|),
\end{split}
\eeq
which can be readily computed for every $\varepsilon$. 
We remark that the second term of \eqref{eq:DKL2} is close to $0$ when the kernel function $k_n(w)$ in \eqref{eq:Galerkin_w} is close to $0$ at every sample $w_m^{l-1}$ for $m \neq n$, so we only need to consider the first term of \eqref{eq:DKL2}. Moreover, to guarantee that $\mathcal{D}_{\text{KL}}^{(2)}$ is reduced for a suitable $\varepsilon$, we can find sample-dependent step sizes $\varepsilon(w_n^{l-1})$ such that 
\beq\label{eq:logposterior_n}
-\log(\pi_y(T_l(w_n^{l-1})))
\eeq
is reduced for each $n = 1, \dots, N$.

\section*{Appendix C: Complexity analysis for parallel pSVN}

We presented a parallel implementation of pSVN in Algorithm \ref{alg:pSVN_parallel}. Lines 4 and 12 involve global communication(gather and broadcast) of the low-dimensional samples $w_m$, $m = 1, \dots, M$, of size $M r$, which are used for the kernel and its gradient evaluations at all samples, as well as for the sample update in \eqref{eq:Galerkin_w}. Line 7 involves global communication (gathers and broadcasts) of the gradients (of size $M r$) and Hessians (of size $M r^2$) of the log posterior density \eqref{eq:gradient_w}, which are used in the expectation evaluation at all samples for assembling the system \eqref{eq:system_w_lumped}. Line 9 involves global communication (gathers and broadcasts) of the kernel values (of size $NM$) at all samples, which are used in moving the samples by \eqref{eq:Galerkin_w}. Meanwhile, Line 9 gathers a local sum of the kernel values $\sum_m k_m(w)$ (of size $N$) and its gradients $\sum_m \nabla_w k_m(w)$ (of size $rN$), performs a global sum of them, and broadcasts the results to all cores, which are used for assembling the lumped Hessian \eqref{eq:Hessian_m}. In summary, the data volumes of communication in Algorithm \ref{alg:pSVN_parallel} are bounded by $\max(M r^2, M N)$ floats.

To implement a parallel version of the adaptive pSVN Algorithm \ref{alg:pSVN_two-level}, we only need to construct the bases $\Psi$ in parallel to replace its Line 5, for which we perform an averaged Hessian action in random directions with $M$ samples in each core by $O(M (r C_h))$ flops, followed by a MPI\_Allreduce with a SUM operator to get a global averaged Hessian action before performing randomized SVD with $O(d r^2)$ flops. The data volumes for communication is $dr$ floats, which dominates all other communication cost if $d$ is so large that $dr > \max(r^2 M, NM)$. Alternatively, we can construct the bases $\Psi$ using Hessian at the local samples in each core without communication for $\Psi$.

\section*{Appendix D: Bayesian Autoencoder Example}

We consider a Bayesian inference problem constrained by a convolutional autoencoder neural network.

In the Bayesian autoencoder problem, we seek to learn a low dimensional representation of data under uncertainty. Given input data $z \in\mathbb{R}^\text{data}$ the $2m$ layer autoencoder mapping is defined as
\beq
y(\cdot) = \circ_{i=1}^{2m} \phi_i(w_i \ast (\cdot) +b_i)
\eeq
where $w_i$ is the convolution kernel (weights) for layer $i$, and $\phi_i$ is an nonlinear activation functions. The $\ast$ operations represents both convolution and downsampling. The first $m$ compositions map down to a low dimensional latent representation of the input data $z$, the last $m$ compositions map the data back to $\mathbb{R}^\text{data}$. 

The data $z$ for the problem are $1000$ randomly selected MNIST images. The target data has 5\% i.i.d. noise added to it based on min-max normalization of the data. The objective function for the autoencoder training problem is a least squares misfit that measures the error between the reconstructed input image and the noisy target image. The inference parameter $\{x_i\} = \{(w_i,b_i)\} \in \mathbb{R}^d$ has the i.i.d. prior $\mathcal{N}(0,\sigma_i^2)$, where $\sigma_1 = 1$, and $\sigma_{i+1}^2 = 0.5\sigma_i^2$. We use a fixed convolution kernel support of $4\times 4$ and vary the number of filters on each layer from $2,4,8$ and use  $m=2$ layers. 

Low rank structure of Hessians has been observed for neural network training problems \cite{AlainRouxManzagol2019,GhorbanKrishnanXiaoi2019,SagunBottouLeCun2016}. Due to the low dimensional nature of the autoencoder, the pSVN algorithm can efficiently find a $r$ dimensional Hessian subspace. The dimensionality of this subspace depends on the decay of the absolute eigenvalues $|\lambda_i|$.

Numerical results are shown below in Figure \ref{AE_figures}. In these trials the problem dimension of the inference parameter is $133$; $128$ particles were used. A fixed candidate rank was chosed to be $r=40$, which is the effective rank of the prior preconditioned Hessian for the problem as seen in the left figure in \ref{AE_figures}. The right figure shows that pSVN minimizes the objective function in training faster than the SVN algorithm for this particular example.

\begin{figure}[!htb]
	\begin{center}
		\includegraphics[width=0.45\columnwidth]{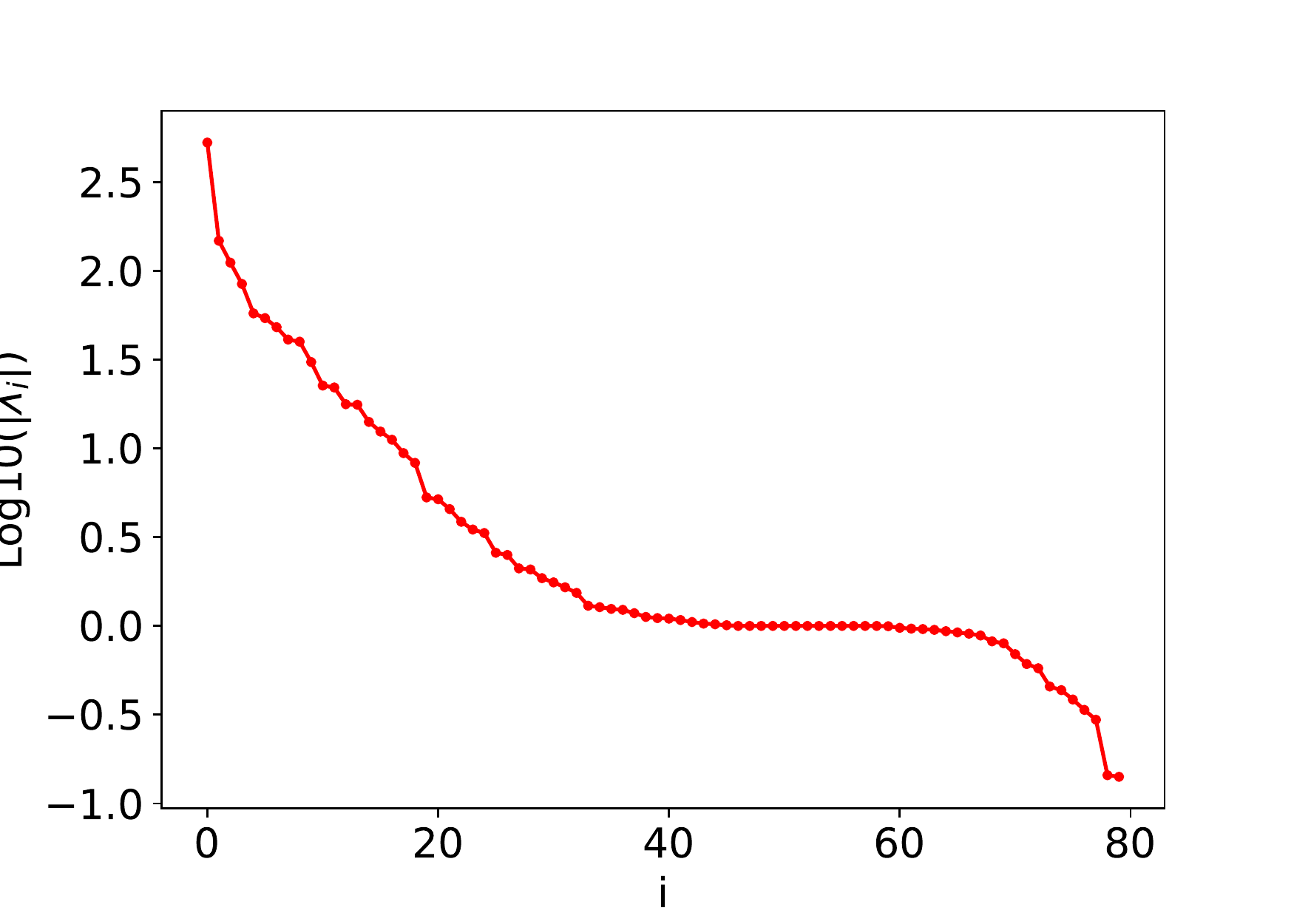}
		\includegraphics[width=0.45\columnwidth]{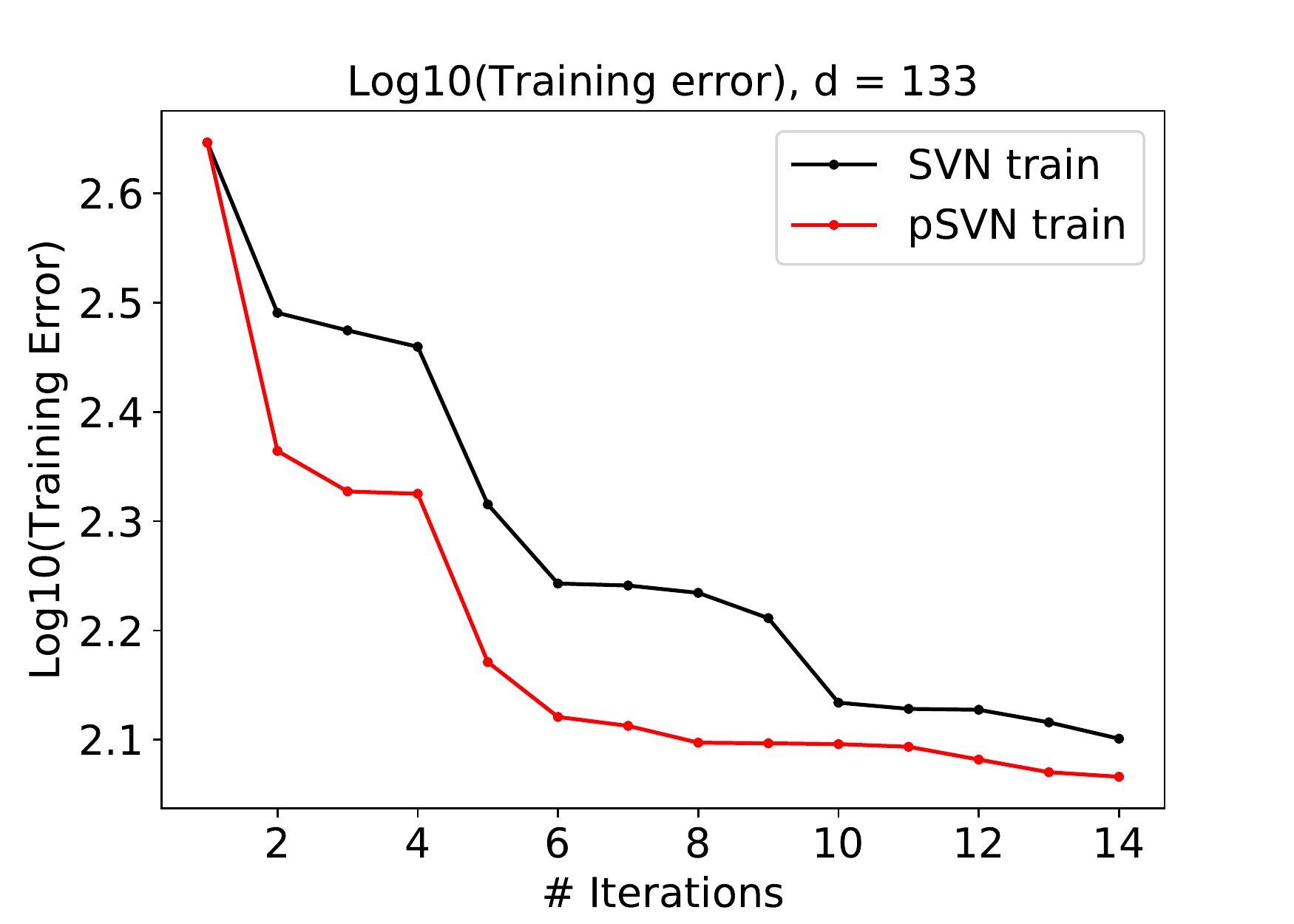}

		\vspace{-0.4cm}
		\caption{Left: Absolute value of eigenvalues of the prior preconditioned Hessian used for the pSVN subspace. Right: Training error for pSVN vs SVN.}
		
		\label{AE_figures}
	\end{center}
	\vspace{-0.4cm}
\end{figure}

\section*{Appendix E: Code}
We implemented the stein variational methods (and the DILI MCMC method) in hIPPYlib (\url{https://hippylib.github.io/}), a python library for solving inverse problems, which relies on FEniCS (\url{https://fenicsproject.org/}), a computing platform for solving partial differential equations. The code for our tests can be downloaded from \url{https://github.com/cpempire/pSVN}.

\end{document}